\input amstex
\documentstyle{amsppt}
\magnification=\magstep0
\define\cc{\Bbb C}
\define\z{\Bbb Z}
\define\r{\Bbb R}

\define\N{\Bbb N}

\define\E{\Cal E}

\define\e{\varepsilon}

\define\CB#1{\Cal C_b(#1)}
\define\st{\subset }
\define\al{\alpha  }
\topmatter
  \title
 Eberlein  almost periodic functions that are not pseudo almost periodic.
  \endtitle
  \author
  Bolis Basit and Hans G\"{u}nzler
\endauthor
 \abstract
 {We construct  Eberlein almost periodic functions $ f_j : J \to H$  so that
    $||f_1(\cdot)||$ is not ergodic and thus not Eberlein almost periodic and $||f_2(.)||$ is Eberlein almost periodic, but  $f_1$ and $f_2$ are not pseudo almost periodic, the Parseval equation for them fails, where $J=\r_+$ or $\r$ and $H$ is a Hilbert space. This answers
several questions posed by Zhang and Liu [18]}
\endabstract
\endtopmatter
\rightheadtext{Eberlein almost periodicity} \leftheadtext{
Basit,   G\"{u}nzler}
  \TagsOnRight
\document
\pageno=1 \baselineskip=20pt

 \head {\S 1.Introduction and Notation}\endhead

Recently Zhang and Liu [18] asked, whether for Hilbert space valued Eberlein
almost periodic $f : \r \to H$ (see \S 2) a Parseval equation holds (Fourier
coefficients for such $f$ are always defined by [14, Theorem 2.4, for $\r_+$]);
this would imply that such $f$ are pseudo almost periodic (see (2.8)).
If additionally the range $f(\r)$ is relatively norm compact, this is true by results
of Goldberg and Irvin [9, Proposition 2.9].

 Here \footnote {AMS subject classification 2010: \,\,\,\, Primary  {43A60}\,\,\,\, Secondary {42A16, 42A75, 42A99}.
 \newline\indent Key words:  Eberlein   almost periodic, pseudo-almost periodic, Parsival equation.}
we show by  examples, that without $f(\r)$ relatively compact the $f$
is in general no longer pseudo almost periodic, one has no Parseval equation.

$\r_+ = [0,\infty)$, $J\in \{\r_+,\r\}$,  $X$ real or complex Banach
space; for  $f : J \to X$,  $f_s(t) : = f(s+t)$, $|f|(t) : = ||f(t)||$, $C_b(J,X)  = \{f :J \to X:$
$f $ continuous, $||f|| =$ sup$_{t\in J}\, ||f(t)||  < \infty\}$, $C_{ub}(J,X) = \{f \in C_b(J,X) : f$ uniformly  continuous$\}$,
$AP(\r,X) = $ almost periodic functions [1, p. 3], [16, p. 18-19], $AP(\r_+,X) =
AP(\r,X)|\r_+ $.

 \head {\S 2.   Eberlein and pseudo almost periodic functions}\endhead

A function $f :J\to X$  is called $Eberlein\,\, almost\,\, periodic$ if $f\in C_b(J,X)$ and orbit  $O(f) : =
\{ f_s : s \in J \}$ is relatively weakly compact in $C_b(J,X)$ ([8, Definition 10.1, p. 232], [6, Definition 1.4], [12, p. 467],  [5, Definition 2.1, p. 138])

  (2.1) \qquad      $ EAP(J,X)  : =  \{ f : f $ Eberlein weakly a. p.$\}$,

 (2.2) \qquad   $EAP_0(J,X) : = \{ f\in EAP(J,X) : 0\in$ weak closure of $O(f) \}$,

 (2.3) \qquad   $EAP_{rc}(J,X) : = \{ f \in EAP(J,X):\,\, $

 \qquad \qquad \qquad  range $f(J)$ relatively norm compact in $X\}$.

 By [2, Theorems 2.3.4 and  2.4.7], [14]  one has

 (2.4)  \qquad     $EAP(J,X)   \st   \E(J,X)  \cap  C_{ub}(J,X)$,  where

 (2.5) \qquad   $ \E(J,X) : = \{f\in L^1_{loc}(J,X) : \text {to\,} f \,\, \text {exists\,\,\,\,}  x\in X$ with

\qquad \qquad \qquad  \qquad
            $||\frac{1}{T}\int _s^{s+T} f(t)\,dt - x|| \to 0$ as $T\to \infty$, uniformly in  $s\in J \}$,

\noindent           then   $m_B(f) : = x$ is called the Bohr-mean.

 For $J$ and $X$ as in \S 1 one has a decomposition theorem [13,
p. 18]
(in $f = g + h$ the $g\in AP(J,X)$, $h\in EAP_0 (J,X)$ are unique)

(2.6) \qquad   $EAP(J,X)   =   AP(J,X)   \oplus   EAP_0 (J,X) $.

The class of $pseudo\,\, almost \,\,periodic$ functions  introduced by Zhang [16],
[17, Definition 5.1, p. 57], [3, (1.1)] is given by

(2.7)\qquad  $PAP_0(\r,X) : = \{f\in C_b(\r,X)$, $m_B(|f|) : = $

\qquad \qquad\qquad\qquad \qquad\,\, $\lim_{T\to \infty} \frac{1}{2T}
                             \int ^T_{-T} |f| (t)dt\,\,\,$  exists $\,\,\, = 0 \}$,

\noindent  similarly for $\r_+$, with $m_B(|f|)=\lim_{T\to \infty} \frac{1}{T}\int ^T_0  |f|(t) dt $,

(2.8)\qquad  $PAP(J,X) : = AP(J,X)  \oplus  PAP_0(J,X) $.

Now by  [9, Proposition 2.9] one has

(2.9)  \qquad
$f\in EAP_{rc}(J,X)$  implies $|f|\in EAP(J,\cc)$,
and  so $|f|^2\in EAP(J,\cc)$

\noindent by [8, Theorem 12.1,  p. 234].

So if $X =$ complex Hilbert space $H$, the polarisation formula ([11, p. 24, (2)]) yields  $(f(\cdot),g(\cdot))_H \in EAP(J,\cc)$ if
$f, g \in EAP_{rc}(J,H)$,  (2.4) shows that  $(f,g) : = m_B (f(.),g(.))_H $ is well defined.
With this (semi-definite) scalar product one gets [9, Theorems 5.2
and 5.7] a Parseval equation for  $f \in EAP_{rc}(J,H)$.

So
(2.6), (2.4) and [9, Corollary 4.19] give

(2.10) \qquad   $EAP_{rc}(\r,H)  = AP(\r,H)  +  \{f\in EAP_{rc}(\r,H) : (f,f) =0\} $,

\noindent with  $(f,f)=0 $ if and only if $ m_B(|f|) = 0$, $f \in  EAP_{rc}(\r,H)$.

So with (2.8) one gets for any complex Hilbert space $H$

(2.11) \qquad  $ EAP_{rc}(\r,H)   \st    PAP(\r,H) $.

Without the "range relatively compact" however all this is no longer true:

\head {\S 3.  Examples}\endhead

For the following we need a converse of Mazur's theorem ([15, p. 120, Theorem 2],
namely

\proclaim{Proposition 3.1} A sequence $(x_n) \st X$  weakly converges to $x\in X$ if and only if, for each subsequence  $(x'_n)$ of  $(x_n)$, there exists a sequence $(y_n)$ of finite convex combinations of the elements of $(x'_n)$  with $||y_n-x||\to 0$ as $n\to \infty$.
\endproclaim
For a proof see [4, Proposition 1.8, p. 17].

\proclaim{Example 3.2} For $J = \r_+ $ or $\r$ and $H$  infinite dimensional
   Hilbert space there exists  $f  \in  EAP_0(J,H)$ so that $m_B(|f|)$ and
  $m_B(|f|^2)$ do not exist, $m_B (|f|)= \lim_{T\to\infty}\frac{1}{T}\int^T_0 ||f(t)||\, dt$ respectively $ \lim_{T\to\infty}\frac{1}{2T}
  \int^T_{-T} ||f(t)||\, dt$ if $J = \r_+ $ respectively   $\r$.  So $|f|$ and  $|f|^2 = (f(\cdot),f(\cdot))_H $ are not in $EAP(J,\cc)$.
\endproclaim
\demo{Proof} Choose an orthonormal sequence $(e_n)_{n \in \N}$ from $H$. Define
$h : \r \to H$  by  $h : = 0$ on $(-\infty, \frac{1}{2}]$, $h(n) : = e_n $, $n \in \N$, $h$ linear on
$[n-\frac{1}{2},n]$ and on $[n,n+\frac{1}{2}]$, with $h(n-\frac{1}{2}) = 0$, $n \in \N$. Then h is well defined
and  $\in C_{ub}(\r,H)$.

Define further $\phi : \r\to [0,1]$  for given $I_n = [\al_n,\beta_n]$, $\beta_n =
\al_{n+1} \in \N$, $\al_n < \beta_n$,  $n \in \N$,  $I_1 = [0,1]$, as follows :

$\phi : = 0$ on
$(-\infty,0]$ and all $I_n$ with odd $n$, $\phi = 1$ on $[\al_{2n} + \frac{1}{10},\beta_{2n}-\frac{1}{10}]$ and $\phi$
linear on $[\al_{2k},\al_{2k} + \frac{1}{10}]$ and $[\beta_{2k}-\frac{1}{10},
\beta_{2k}]$, $k \in \N $. Then also $\phi$ is well defined and  $\in C_{ub}(\r,\r)$.

To get a non-ergodic $\phi$, choose the $I_n$ recursively with $I_1= [0,1]$ as follows
(Zorn's Lemma):

 If $I_1,\cdots ,I_{2k}$ are defined, take $\al_{2k+1}: = \beta_{2k}$
and $\beta_{2k+1}$ such that  $\frac{\al_{2k} }{ \beta_{2k+1}} < \frac{1} {5}$;

if $I_1, \cdots, I_{2k-1}$
are defined, take $\al_{2k} : = \beta_{2k-1}$ and $\beta_{2k}$ such that

$\frac{\beta_{2k} - \al_{2k} - \frac{1}{5}  - 2}{  \beta_{2k} } > \frac{3}{4} $.

Finally, define   $f  : =  \phi h| J$,  $\in  C_{ub}(J,H)$. Then

$\lim\, \inf_{T\to \infty} \frac{1}{2T} \int^T_{-T} |f| (t)\, dt = lim\, \inf _{T\to \infty}\frac{1}{2T} \int^T_0 |f|(t)\,dt\le $

     $\lim\, \inf_{T\to \infty} \frac{1}{2T} \int ^T_0 \phi(t) \, dt \le \lim \, \inf_{k\to \infty}\frac {1}{\beta_{2k+1}}
       \int_0^{\beta_{2k+1}} \phi (t)\, dt \le \frac {1}{5}$,

  $\lim\, \sup_{T\to \infty} \frac{1}{2T} \int^T_{-T} |f|^2(t) dt = lim sup_{T\to \infty} \frac{1}{2T} \int^T_0
      ||\phi(t) h(t)||^2\, dt \ge$

      $  \lim \,\sup_{k \to\infty}\frac {1}{\beta_{2k}} \int_0^{\beta_{2k}} |\phi|^2(t)
      |h|^2 (t)\, dt \ge $

      $ \lim\, \sup_{k\to \infty}\frac {1}{\beta_{2k}} [\beta_{2k} - \al_{2k} -\frac{1}{5} -2] 1^2 \frac{1}{3}
      \ge    \frac{1}{4} $.

Since  $|f|^2  \le |f|$, the above shows that  $m_B( |f| )$  and  and $m_B(|f|^2)$
   do not exist, $J =\r$ or $\r_+ $.

 $f  \in EAP(J,H)$ : With the Eberlein-Smulian theorem [7, p. 430, 1.
Theorem] one has to show : to each sequence $(b_n)$ from $J$ there
exists a subsequence $(a_m)$ and $g \in  V : = C_b(J,H)$ with $f_{a_m}\to g$
  weakly in V. Now if $(b_n)$ is bounded, there exists a subsequence $(c_k)$
 and $c \in J$ with $c_k\to c$, so  $f_{c_k}\to f_c$ uniformly on $J$  and so even in
the norm of V, since  $f \in  C_{ub}(J,H)$.

 Now assume  $a_m\to \infty $; by taking a further subsequence, one can assume
 $a_{m+1} - a_m > 1$ and $a_1 \ge 1$, $m \in \N$. To apply Proposition  3.1,  let $(c_k)$ be any subsequence of $(a_m)$ and $\e  > 0$, then
there exist  $q$, $k_1,..,k_q  \in \N$  with $\frac{1}{q} < \e^2$, $c_{k_{j+1}}  - c_{k_j}
>  1 $ and $c_{k_1} \ge  1$, $1 \le j \le q$. Then the $c_{k_j}$ are in different
$[n-\frac{1}{2},n+\frac{1}{2}]$ intervals for different $j$, so for any $t \in \r$,

   $f_{c_{k_j}}(t) = f(c_{k_j} + t) = r_{j,t} e_{p(j,t)}$  with
   $0\le r_{j,t}\le 1$ and  $p(i,t)  < p(j,t)$  if  $i  < j$  and

    $c_{k_ i} + t > \frac{1}{2} $.
With $i_0$ minimal if such $i$ exist and  $\theta_{k_j} = \frac{1}{q}$,  $1\le j\le q$, else $= 0$,

\noindent one gets

  $||\sum_{j=1}^{k_q} \theta_{k_j} f(c_{k_j} + t)||_H ^2  = ||\sum^q_{j=1} \frac{1}{q}
      f(c_{k_j} + t) ||_H ^2  = ||\sum^q_{j=i_0} \frac{1}{q} r_{j,t} e_{p(j,t)}||_H^2 =$

      $
       \frac{1}{q^2} \sum^q_{j=i_0} (r_{j,t})^2 $
      $\le $
      $\frac{q}{q^2}  = \frac{1}{q}$;

\noindent  if no such $i_0$ exists, the above sum is even $0 \le \frac{1}{q} $.

  This holds for any  $t \in\r$, so  $||\sum _{j=1} ^{k_q} \theta _{k_j} f_{c_{k_j}}||_V  < \e$.
  Therefore by Proposition 3.1 indeed  $f_{c_k}\to 0 $ weakly in $V$,
   $J =\r$ or $\r_+$.

  The case  $a_m\to - \infty$ ( $ J =\r$) follows similarly.

$ f  \in  EAP_0(J,H) $: For $(b_n) = (n)$ the above shows  $0 \in $ weak closure
     of orbit $O(f)$.

  $|f|$ and $(f(.),f(.))$ not Eberlein almost periodic follows with (2.4).  \P
\enddemo

Since for the $f$ of Example 3.2 the Bohr mean $m_B(|f|^2)$ does not exist, one has
no Parseval equation.

$ EAP(J,X) \st PAP(J,X)$ is also false, already for X = Hilbert space :

Assume $f\in EAP(J,X) \st PAP(J,X)$. Then  $ f = g+h$, $g \in AP(J,X)$, $h \in PAP_0(J,X)$;
now for  $f \in EAP_0(J,X)$ one can show that all Fourier coefficients vanish (for J = $\r_+ $ see
[14, Theorem 2.4]), for $h$ the same holds , implying $g=0$, then $f=h \in PAP_0(J,X) $ and so the
existence of $m_B(|f|)$, a contradiction for  $f$ of Example 3.2.

 The proof of Example 3.2 works also for $X = l^p(N,\cc)$, $1<p \le \infty$ and $c_0$, so
$EAP (J,X) \st  PAP(J,X)$ is also false for these $X$.

Since for any $f  \in  EAP(J,X)$ the range $f(J)$ is relatively weakly compact,
  and if $X = l^1 = l^1(M,\cc)$, any $M$, this implies $f(J)$ relatively norm compact [10, p. 281 (2)],  one
has

 $EAP (J,l^1)=EAP_{rc} (J,l^1) \st  PAP(J,l^1)$.

\proclaim {Example 3.3}  For $J =\r$ or $\r_+$ and $H$  separable infinite dimensional Hilbert space there
   exist  $ f  \in  EAP_0(J,H)$ with  $|f|, |f|^2   \in  AP(J,\r) \st  EAP(J,\r)$, but  $f(J)$ is not
relatively compact, $m_B(|f|)$ and $m_B(|f|^2)$ exist  and are $ > 0$.

So a converse of (2.9) is not true, even with $|f|, |f|^2  \in EAP (J,\r)$ the Parseval equation can fail, such f need not be pseudo almost periodic.
\endproclaim

\demo{Proof}  Choose an orthonormal sequence $(e_n)_{n \in \z}$ from $H$. Define
$f : \r \to H$  by   $f(n) : = e_n $, $n \in \z$, $f$ linear on
$[n-\frac{1}{2},n]$ and on $[n,n+\frac{1}{2}]$, with $f(n-\frac{1}{2}) = 0$, $n \in \z$. Then f is well defined
and  $\in C_{ub}(\r,H)$. One can prove that $f \in EAP_0(\r,H)$ as in the
proof of Example 3.2. Obviously, $|f|  \in  C_{ub}(\r,\r)$ has period  $1$ and so  $|f| \in AP(\r,\r) \st EAP(r,\r)$.
$f|\r_+ $ has the same properties.   \P
\enddemo

\bigskip

\indent School of Math. Sci., Building 28M,
 Monash University, Vic. 3800,
 Australia.

\indent E-mail: bolis.basit\@sci.monash.edu.au.

\bigskip

\indent Math. Seminar der  Univ. Kiel, Ludewig-Meyn-Str.,
 24098 Kiel, Deutschland.

\indent E-mail: guenzler\@math.uni-kiel.de.


\Refs

\ref\no1\by L. Amerio and G. Prouse \book Almost-Periodic
Functions and Functional Equations, Van Nostrand, 1971
\endref
\ref\no2\by B. Basit\book Some problems concerning different types
of vector valued almost periodic functions,  Dissertationes Math.
338 (1995), 26 pages
\endref
\ref\no3\by B. Basit and C. Zhang\book  New almost periodic type
functions and solutions of differential equations,
 Can.J.Math. 48 (1996), 1138-1153
\endref
\ref\no4\by B. Basit,  and  A.J. Pryde\book Asymptotic behavior of $C_0$-semigroups and solutions of linear  and semilinear abstract differential equations, Russ. J.  Math. Physics, 13 (1), (2006), p. 13-30
\endref
\ref\no5\by J. F. Berglund, H. D. Junghann and P. Milnes\book
Analysis on Semigroups, John Wiley,  1989
\endref
\ref\no6\by Burckel\book Weak almost periodic functions on semigroups, Gordon and Breach Science Publishers, New York, 1970
\endref
\ref\no7\by N. Dunford and J. T. Schwartz   \book Linear
Operators, Part I, Interscience, New York, 1963
\endref
\ref\no8\by W. F. Eberlein\book Abstract  ergodic theorems and weak almost periodic functions, Trns. Amer. Math. Soc., 67 (1949), 217-240
\endref
 \ref\no9\by  S. Goldberg and P. Irwin   \book Weakly almost
periodic vector-valued  functions,  Dissertationes Math. 157
(1979), 42 pages
\endref
\ref\no10\by G. K\"{o}the\book Topological Vector Spaces I, Springer 1969
\endref
\ref\no11\by L. H. Loomis\book An Introduction to  Abstract Harmonic Analysis, D.  Van Nostrand Co., New York, 1953
\endref
\ref\no12\by P. Milnes\book On vector-valued weakly  almost
periodic functions, J. London  Math.  Soc. (2) 22 (1980), 467-472
\endref
\ref\no13\by W. M. Ruess and W. H. Summers\book Integration of
asymptotically almost periodic  functions and
 weak almost periodicity, Dissertationes Math. 279  (1989), 35 pages
\endref
\ref\no14\by W. M. Ruess and W. H. Summers\book Ergodic theorems
for semigroups of operators,
 Proc. Amer. Math. Soc. 114  (1992), 423-432
\endref
\ref\no15\by K. Yosida \book Functional Analysis, Springer Verlag,
  1976
\endref
\ref\no16\by C. Zhang\book  Vector-valued means and their applications in some
vector-valued function spaces,
 Dissertationes Math. 334    (1994), 35 pages
\endref
\ref\no17\by C. Zhang\book  Almost Periodic Type  Functions and
Ergodicity,
 Science Press/Kluwer Acad. Publ., 2003
 \endref
\ref\no18\by C. Zhang and  J. Liu\book  Two unsolved problems on almost periodic type functions, Appl. Math. Letters 23 (2010), 1133 - 1136
 \endref
\endRefs
\enddocument